\documentclass[letterpaper, 10 pt, conference]{ieeeconf}  

\IEEEoverridecommandlockouts
\usepackage[dvipdfm]{graphicx} 
\usepackage{subfigure}
\usepackage{amsmath}
\usepackage{amsfonts}
\usepackage[psamsfonts]{amssymb}
\usepackage{color}
\usepackage{comment}
\usepackage{times} 
\usepackage{dsfont}
\usepackage{graphicx}
\usepackage{epsfig}
\usepackage{epstopdf}

\title{\LARGE \bf
Temperature Regulation in Multicore Processors Using  Adjustable-Gain Integral Controllers
}

\author{K. Rao, W. Song, S. Yalamanchili, and Y. Wardi$^*$$^{\dag}$
\thanks{$^*$School of Electrical and Computer Engineering, Georgia Institute of Technology,
Atlanta, GA 30332. Email: raokart@gatech.edu, wjhson@gatech.edu, sudha@ece.gatech.edu, ywardi@ece.gatech.edu.}\thanks{$^{\dag}$Research supported in part by  NSF under Grant Number CNS-1239225.}
}

\begin{document}

\maketitle

\begin{abstract}
 This paper considers the problem of temperature regulation in multicore processors by dynamic voltage-frequency scaling. We propose a  feedback
  law that is based on an integral controller with adjustable gain, designed
  for fast tracking convergence  in
  the face of model uncertainties, time-varying plants, and tight
  computing-timing constraints. Moreover, unlike prior works we
  consider a nonlinear, time-varying plant model that trades off precision for simple and efficient on-line
  computations. Cycle-level, full system simulator implementation and
  evaluation illustrates fast and accurate tracking of given
  temperature reference values,  and compares favorably with fixed-gain
  controllers.

\end{abstract}

\section{Introduction}
\label{Section1}
The end of Dennard scaling has led to increasing power densities on
the processor die and consequently higher chip temperatures~\cite{dennard1974design, ITRS2011}.
Emerging and future processors are thermally limited and must operate
within the cooling capacity of the chip package, which is typically
represented by the maximum operating temperature. Dynamic Thermal
Management (DTM) techniques have emerged to manage thermal behaviors
and are challenged by a number of phenomena. In particular, the exponential dependence of static power on temperature limits the effectiveness of many existing DTM techniques.  This
coupling can also lead to thermal runaway that must be prevented by
DTM to avoid damaging the chip. Furthermore, the structure of the
thermal field matters as spatial and temporal variations in the
thermal field degrade device reliability and accelerate chip
failures. Similarly, rapid changes in the thermal field referred to
as thermal cycling, also cause thermal stresses that degrade device
and hence chip reliability.

A specific class of thermal regulation techniques includes   activities' management  like   instruction fetch
throttling and clock gating ~\cite{skadron2003temperature, skadron2002control}, thread migration (computations' rescheduling) ~\cite{liu2012neighbor, yeo2008predictive}, and core frequency scaling
~\cite{kong2012recent}. References ~\cite{skadron2003temperature, skadron2002control}
use PI and PID controls to slow down the rate of the instruction-fetch unit whenever the temperature exceeds a given upper bound,
while ~\cite{liu2012neighbor, yeo2008predictive} schedule threads (computations) from hot cores to cooler cores in effort to
maintain a balanced thermal field.
Initial heuristic approaches  started giving way to control-theoretic formalisms, with  the aforementioned references
~\cite{skadron2003temperature, skadron2002control} providing (to our knowledge) the earliest examples. Subsequently,
Reference
\cite{donald2006techniques} considered a similar upper-bound regulation problem but uses Dynamic Voltage Frequency
Scaling (DVFS) for temperature control.
More recently~\cite{qian2011cyber} described
a controller for regulating the fluid in a microfluidic heat sink
based on the measured temperature as well as predicted temperature
estimated from the projected power profile. Other work has investigated DTM
under soft and hard real-time constraints~\cite{shi2010dynamic, fu2012feedback} seeking to satisfy
thermal upper bounds  while operating under scheduling constraints.

More recently, there emerged a number of approaches, which are based
on optimal control and optimization. Reference~\cite{zanini2009multicore}
minimizes a least-square difference between the working frequency and the
frequency mandated by the operating system, subject to thermal and frequency
constraints, by using model-predictive control. Reference~\cite{wang2011adaptive}
uses similar techniques to minimize the least-square difference between set power
levels and actual power levels in a core. Reference~\cite{murali2008temperature}
uses a combination of off-line convex optimization and on-line control  to obtain
uniform spatial temperature gradient across several cores in a processor. We point out
that these references assume linear and time-invariant plant-models for their
respective control systems;~\cite{wang2011adaptive}
updates the model on-line while~\cite{zanini2009multicore,murali2008temperature}  do not.
Finally, reference~\cite{bartolini2013thermal} minimizes energy consumption while preserving
performance levels within a tolerable limit by employing separate Model Predictive Controllers
for each core to ensure thermal safety, and updates the power-temperature model for the cores online.

Besides the need to limit  core and chip temperatures, there
 is a  pressure to maintain temperatures close to
package capacity in order to maintain high levels of performance. This typically is achieved by  adjusting
the rates of the processor cores as, for example, in
Intel processors ~\cite{rotem2012power} and  AMD
processors ~\cite{AMDPhenom}. Moreover,
spatiotemporal variations in the thermal field generally impact device
degradation and energy efficiency. For example, thermal gradients
between adjacent cores on a die increase leakage power in the cooler
core, thereby increasing its temperature and reducing its energy efficiency
(ops/joule)~\cite{paul2013}. Further, the stresses introduced by the
gradients reduce lifetime reliability by accelerating device
degradation~\cite{song2015}. These affects are exacerbated in
heterogeneous multicore processors where cores of different
complexities (and therefore thermal properties) are utilized to
improve overall energy efficiency. Consequently, it has become
necessary to be able to allocate and control the usage of thermal
capacity in different regions of the die. Core-temperature regulation (and not only optimization)
can provide an important means to this end.

This paper proposes an approach for regulating core temperatures by DVFS so as to track
given reference temperature values (set points). The
frequency is adjusted by an integral controller with adjustable gain,
designed for fast tracking-convergence under changing program
loads. Unlike the aforementioned references that are based on optimal control and optimization,
we consider a
nonlinear, time-varying plant model that captures the exponential dependence of temperature on static
power. The basic idea is to  have
the on-line computations of the integrator's gain be as simple and efficient as possible even at the expense of precision.
This is made possible by a great degree of  robustness of the tracking performance of the controller with
respect to variations from  the designed
integrator's gain, which was observed from extensive simulations (see \cite{almoosa2012power} for
analysis and discussion). We verify the efficacy of our technique by  simulations on   a full
system, cycle level simulator executing industry standard benchmark
programs, and demonstrate rapid convergence
despite the modeling errors and changing program loads.

We first applied  the proposed approach in   \cite{almoosa2012power}
for controlling the dynamic core power via DVFS. The problem considered here is more challenging
for the following two reasons. 1).  The underlying model required in this paper is much
more complicated. Ignoring the static power permitted Reference \cite{almoosa2012power} to use an
established third-order polynomial formula for the dynamic power as a function of frequency. In contrast,
the temperature's dependence on frequency has no explicit formula, but  rather is described implicitly by a
differential equation that models the heat flow. Furthermore, the temperature depends on the total (static and dynamic)
power while the static power depends on the temperature (and voltage), and this circular dependence was avoided  in \cite{almoosa2012power} by ignoring the static power.\footnote{In present-day  technologies and applications
the static power can be as high as the dynamic power and  no-longer can be ignored.} For reasons discussed later,
the duration of the control cycle is about $10$ms, which requires fast computations in the loop. Our main challenge in this regard was to find an approximate model yielding simple computations while preserving the aforementioned convergence properties of the
control algorithm.
2). The temperature levels in different cores on a chip are inter-related due to the diffusion of heat between them,
while their dissipated dynamic powers are not directly related to each other by such physical laws. Therefore it is natural for the
dynamic-power control law in \cite{almoosa2012power} to be distributed among the cores, while in this paper the temperature control appears to have to be centralized. Nonetheless we argue for a distributed control law and justify its use via analysis and simulation.

The next section presents our regulation techniques  in an abstract setting  and recounts relevant existing results. Section III describes  our modeling approach to the thermal regulation problem, Section IV presents simulation results
on standard industry benchmarks, and Section V concludes the paper.

\section{Regulation Technique}

Consider the discrete-time, Single-Input-Single-Output (SISO)
feedback system shown in Figure 1,
whose input is a constant  reference $r$,  its output is denoted by $y_{n}$,
the input to its controller is the error signal $e_{n}$, and the input to the plant is $u_{n}\in \mathds{R}$.
Suppose  that the plant is a  time-varying  nonlinear system  described via the relation
\begin{equation}
\label{y_g_u}
y_{n}=g_{n}(u_{n-1}),
\end{equation}
where the function $g_{n}:\mathds{R}\rightarrow \mathds{R}$ is called the {\it plant function}.

\begin{figure}[!t]
\centering
\includegraphics[width=3in]{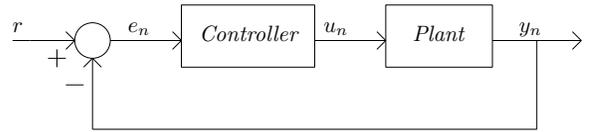}
\caption{Control System  Block Diagram}
\end{figure}

If the controller is an integrator having the transfer
function $G_{c}(z)=Az^{-1}/(1-z^{-1})$, for a constant $A>0$, then in the time domain it is defined by  the relation $u_{n}=u_{n-1}+Ae_{n-1}$. However, we will consider an adjustable (controlled)
gain, and hence the controller equation has the form
\begin{equation}
\label{ControlEquation}
u_{n}=u_{n-1}+A_{n}e_{n-1},
\end{equation}
where the gain $A_{n}$ is computed in a manner described below. The error signal has the form
\begin{equation}
\label{ErrorEquation}
e_{n}=r-y_{n}.
\end{equation}
Suppose that the plant functions $g_{n}(u)$ are differentiable, and let ``prime'' denote their derivatives with respect to
$u$.
We define the gain  $A_{n}$ as
\begin{equation}
\label{ControllerGain}
A_{n}=\frac{1}{g_{n}^{\prime}(u_{n-1})}.
\end{equation}

The systems considered in the sequel have the following structure. Consider a SISO dynamical system having
an input $\{u(t)\}$ and output $\{y(t)\}$, $t\geq 0$. Partition the time-horizon $\{t\geq 0\}$ into
consecutive time-slots $[\tau_{n-1},\tau_{n})$, $n=1,2,\ldots$, with $\tau_{0}:=0$ and $\tau_{n+1}>\tau_{n}$
$\forall\ n=1,\ldots$; define $C_{n}:=[\tau_{n-1},\tau_{n})$ and call it the $n^{th}$ control cycle. Suppose that
the value of the input is changed only at the boundary points $\tau_{n}$, and denote the value of the input
$u(t)$ during $C_{n}$ by $u_{n-1}$. Let $y_{n}$ be a quantity of interest that is generated by the system
during $C_{n}$ from $u_{n-1}$, such as $y(\tau_{n}^{-})$ or
$\int_{C_{n}}y(t)dt$. $y_{n}$ also depends on the initial condition $y(\tau_{n-1})$, but this is reflected in Equation (1)
by the system's definition as time varying. Thus, (1) represents certain input-output properties of dynamical
systems while hiding the details of the dynamics and  appearing to have the form of a memoryless nonlinearity.
Regarding the feedback system, we suppose that  $u_{n-1}$, $y_{n-1}$, and $e_{n-1}$ are
available to it at time $\tau_{n-1}$, and  it generates  $y_{n}$ by (1) and computes $A_{n}$ during $C_{n}$ via (4).
The  closed-loop system is defined by repeated applications of  Equations $(1)\rightarrow(4)\rightarrow(2)\rightarrow(3)$.

To see the rationale behind the definition of the gain $A_{n}$ in (4) consider the case where
the plant is time invariant, namely  $g_{n}(u)=g(u)$ for a function $g:\mathds{R}\rightarrow \mathds{R}$.
Then this control law amounts to a realization of the Newton-Raphson method for solving the equation
$g(u)=r$, whose convergence means that $\lim_{n\rightarrow\infty}e_{n}=0$.
Furthermore, if the derivative $g^{\prime}(u_{n-1})$ cannot be computed exactly,   convergence  also
is ensured under broad assumptions. For instance, suppose that Equation (\ref{ControllerGain}) is replaced by
\begin{equation}
A_{n}=\frac{1}{g^{\prime}(u_{n-1})+\xi_{n-1}},
\end{equation}
where the  error term $\xi_{n-1}$ is due to modeling uncertainties, noise, or computational errors. If the function
$g(u)$ is globally monotone increasing or monotone decreasing,  and  convex or  concave throughout $\mathds{R}$,
and if the relative error term $|\xi_{n}|/|g^{\prime}(u_{n})|$ is upper-bounded by a constant $\alpha\in(0,1)$
 for all $n=1,2,\ldots$, then convergence (in the sense that $\lim_{n\rightarrow\infty}e_{n}=0$) is
guaranteed for every starting point $e_{0}$ as long as $g^{-1}(r)\neq\emptyset$. If $g(u)$ is piecewise monotone and piecewise convex/concave then convergence is guaranteed for a local domain of attraction; namely, for every
point $\hat{u}\in \mathds{R}$ such that $g(\hat{u})=r$ and $g^{\prime}(\hat{u})\neq 0$, there exists an open interval $I$ containing $\hat{u}$ such that,
for every $u_{0}\in I$, $u_{n}\rightarrow\hat{u}$ and hence $e_{n}\rightarrow 0$ as $n\rightarrow\infty$. More specifically,
there exist $\gamma\in(0,1)$ and $N\geq 0$ such that, for every $n\geq N$,
\begin{equation}
|e_{n}|\leq\gamma|e_{n-1}|.
\end{equation}
These, and more extensive results
concerning convergence of Newton-Raphson method for finding the zeros of a function can be found in \cite{lancaster1966error}.

In the general time-varying case where  the plant function $g_{n}$ is $n$-dependent (as in (1)),  it cannot be
expected to have
$e_{n}\rightarrow 0$. However, the term
$\limsup_{n\rightarrow\infty}|e_{n}|$ has been shown to be bounded by quantified measures of the system's time-variability. For instance, \cite{almoosa2012power} derived the following result under conditions of monotonicity and strict convexity of the functions
$g_{n}$:
For every $\varepsilon>0$ there exist  $\delta>0$
such that, if $|g_{n-1}(u_{n-1})-g_{n}(u_{n-1})|<\delta$ $\forall n=1,2,\ldots$, then
$\limsup_{n\rightarrow\infty}|e_{n}|<\varepsilon.$ Moreover, there exist $\eta>0$ and $N\geq 0$ such that,
for every $n\geq N$, Equation (6) holds true as long as $|e_{n-1}|>\eta$.

These results have had extensions to the multivariable
case arising in Multi-Input-Multi-Output (MIMO)  systems with the same number of outputs as inputs (e.g., \cite{lancaster1966error,ortega2000iterative}).
Accordingly, for a given $M\geq 1$, let $u\in \mathds{R}^M$ and $y\in \mathds{R}^M$ denote the input and output of the plant, respectively. Define the plant function by Equation (1)  except that $g_{n}$ is a function from
$\mathds{R}^M$ to $\mathds{R}^M$, the feedback equation by (2)  except that $A_{n}$ is an  $M\times M$ matrix, the error term
via Equation (3), and the gain matrix $A_{n}$ by the following extension of Equation (4),
\begin{equation}
A_{n}=\Big(\frac{\partial g_{n}}{\partial u}(u_{n-1})\Big)^{-1}.
\end{equation}
In the time-invariant case where $g:=g_{n}$ is independent of $n$, the system consisting
 of  repetitive applications of Equations
$(1)\rightarrow(7)\rightarrow(2)\rightarrow(3)$ comprises an implementation of Newton-Raphson method for solving the equation
$g(u)=r$.

We are concerned with the time-varying case where the plant function depends on $n$ as in (1),
and the Jacobian  matrix $\frac{\partial g_{n}}{\partial u}(u_{n-1})$ is approximated rather than computed
exactly. In this case  Equation (7) is replaced by the following extension of (5),
\begin{equation}
A_{n}=\Big(\frac{\partial g_{n}}{\partial u}(u_{n-1})+\xi_{n-1}\Big)^{-1},
\end{equation}
where the error term $\xi_{n-1}$ is an $M\times M$ matrix.  Define
the relative error at the $n^{th}$ step of the control algorithm by
${\cal E}_{n}:=||\xi_{n-1}||\big(||\frac{\partial g_{n}}{\partial u}(u_{n-1})||\big)^{-1}$.
Various general results concerning  the Newton-Raphson  method guarantee local convergence of the control algorithm
under the condition that
${\cal E}_{n}\leq\alpha$ for some $\alpha<1$, for all $n=1,2,\ldots$; see, e.g.,  \cite{lancaster1966error}.
They typically state that
$\lim_{n\rightarrow\infty}e_{n}=0$ in the time-invariant case, and show upper bounds on $\lim\sup_{n\rightarrow \infty}||e_{n}||$ in
the case of time-varying systems.

The control law defined by Equations (8) and (2) updates all of the $M$ components of $u_{n}$ simultaneously and hence
can be viewed as centralized. However, by ignoring the off-diagonal terms  of
$\frac{\partial g_{n}}{\partial u}(u_{n-1})$ we effectively obtain a distributed controller. Formally, define $D_{n}$ to be the matrix comprised of the diagonal
elements of $\frac{\partial g_{n}}{\partial u}(u_{n-1})$, and define
$\xi_{n-1}:=D_{n}-\frac{\partial g_{n}}{\partial u}(u_{n-1})$. Then Equation (8) can be computed
in parallel by Equation (5) for each input-output coordinate. Thus the system comprised of repeated
applications of Equations $(1)\rightarrow(8)\rightarrow(2)\rightarrow(3)$ can be viewed as a distributed system consisting of
repeated runs of $(1)\rightarrow(5)\rightarrow(2)\rightarrow(3)$.

\section{Temperature Control in Multi-Core Computer Processors}

This section describes an application of the control technique described in Section II to
temperature regulation in computer cores by adjusting their  frequencies.
Unlike the case of regulating the dynamic power, described in \cite{almoosa2012power},
the frequency-to-temperature relationships are highly dynamic and complex, and moreover, the temperatures at various cores on a
chip are inter-related.  Nevertheless  our objective is to have a distributed controller whose required calculations
are as simple as possible
 since, among other reasons, their complexity  poses a lower bound on the durations of the control
cycles.

To this end  we consider approximations that trade off precision with
low computational complexity by leveraging  the  convergence robustness reflected in Equations (5) and
(6). Therefore  much of the developments in this section concern modeling approximations that yield simple computations. The resultant control law is tested
in the next section.

The first part of the investigation concerns the frequency-to-temperature relations in a single core, formalized via  the scalar-version of Equation (1). Suppose that the frequency applied to the core has a constant
value during each control
cycle and it is changed only at the cycle boundaries.
Let $\phi$ denote the frequency applied to the core during a typical control cycle, and let $P:=P(t)$ and $T:=T(t)$ denote the
resulting dissipated power and spatial average temperature during the cycle.
The power has two main components: static power and dynamic power, respectively denoted by
$P_{s}$ and $P_{d}$. The static power is  dissipated due to leakage currents in the transistors, and the dynamic power
 is dissipated when the transistors are switched between the {\it on} and {\it off} states.
 Figure 2 depicts the functional relations between these quantities, and we note that the dynamic power depends on
 the frequency, the temperature depends on the total power, and   the static power depends on the  frequency and temperature.
 The relationships between these quantities are indicated in the figure by the system-notation $S_{1}$, $S_{2}$,
 and $S_{3}$, and we next describe their models in detail.

\begin{figure}
\centering
\includegraphics[width=3.2in]{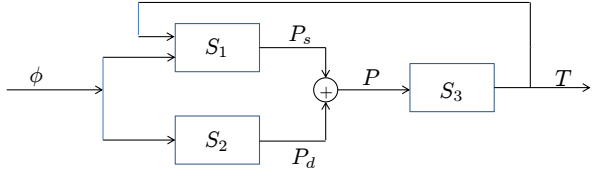}
\caption{System Model}
\label{Figure 2}
\end{figure}

 The core frequency typically is controlled by an applied voltage $V$,
 not shown in Fig. \ref{Figure 2}. The relationship between frequency and voltage
 can be modeled by the affine equation
 \begin{equation}
 \label{VoltageEquation}
 V=m\phi+V_{0},
 \end{equation}
 \cite{burd2000dynamic, mcgowen2006power}
 whose slope $m$ often can be obtained from the manufacturer.

 As mentioned earlier, the total power is given by
 \begin{equation}
 \label{PowerEquation}
 P\ =\ P_{s}+P_{d}.
 \end{equation}
 {\it The system $S_{1}$ (Figure 2):}
 An established physical model for the static power is described in~\cite{butts2000static}, and it is given by the
 equation
 \begin{eqnarray}
 \label{StaticPowerEquation1}
 P_{s}\ =\
 VNk_{{\rm design}}I^{\prime}_{so}e^{-(V_{{\rm off}})q/(\eta kT)}\nonumber \\
 \times 10^{-(V_{T})q/(2.303\eta kT)},
 \end{eqnarray}
 where $V$ is the applied voltage, $N$ is the number of transistors in the core,
 $k_{{\rm design}}$ is a positive parameter depending on the core design,
 $I^{\prime}_{so}$ is a constant related to the subthreshold drain current,
 $V_{{\rm off}}$ is an empirically determined model parameter, $q=1.6\times 10^{-19}$C is the electron's charge,
 $\eta$ is a technology-dependent parameter,
 $k=1.38\times 10^{-23}m^2kgs^{-2}K^{-1}$ is the
 Bolzmann's constant, $T$ is the core temperature in
 Kelvin, and $V_{T}$ is the threshold voltage of the transistor.
 Grouping terms and defining
 \[
 \beta=Nk_{{\rm design}}I^{\prime}_{so}
 \]
 and
 \[
 \gamma=q(V_{{\rm off}} + V_{T})/(2.303\eta k),
 \]
 we obtain the equation
 \begin{equation}
 \label{StaticPowerEquation2}
 P_{s}\ =\ V\beta\times 10^{-\gamma/T},
 \end{equation}
 where we note that $\beta>0$ and $\gamma>0$.
 Observe that $P_{s}$ depends on $V$ (and hence on
 $\phi$ via (9)) as well as on $T$.

 {\it The system $S_{2}$:} An established model for the dynamic power~\cite{rabaey2009low}
 is described by the following equation,
 \begin{equation}
 \label{DynamicPowerEquation}
 P_{d}\ =\ \alpha(t)CV^{2}\phi,
 \end{equation}
 where $C$ is the lumped capacitance of the core, and $\alpha(t)$, called the
 {\it activity factor}, is a time-varying parameter related to the amount of switching activity of the logic gates at the
 core. We note that $\alpha(t)$ cannot be effectively computed or predicted in real time, but its evaluation is not needed
 for the control algorithm.

 {\it The system $S_{3}$:} A detailed physical model for the power-to-temperature relationship
 is quite complicated. However, it will be seen that what we need is
 the derivative term $\frac{dT}{dP}$, and that this can be approximated
 by a constant which can be computed off line. In making this  approximation we leverage
   the robustness of the tracking algorithm
 with respect to errors in the computation of
 $g_{n}^{\prime}(u_{n-1})$ (see (5),(6)), as discussed in Section II.

 The power-to-temperature  relationship in a core  has had an effective model in~\cite{han2007tilts}, that is based on a linear and time-invariant system, and hence yields fast simulation-response time as compared to physics-based models. The dimension of the system is the number of
 functional units in the core, typically in the $50$ - $100$ range,
 the input $u$ represents the vector of the dissipated  power at each functional unit, and the
 state variable $x$ is the temperature at each functional unit. The  state equation has the
 form
 \begin{equation}
 \label{TILTSEquation}
 \dot{x}\ =\ Ax+Bu,
 \end{equation}
 where the matrices $A$ and $B$ can be estimated off line. At each time $t$,  the total dissipated power at the core,
 $P:=P(t)$, and the spatial average of the core temperature,
 $T:=T(t)$, are
 linear combinations of $u$ and $x$, respectively, and therefore the $P-T$ relationship can be described via the scalar
 differential equation
 \begin{equation}
 \dot{T}=aT+bP.
 \end{equation}
 Consequently, the derivative term $\frac{dT}{dP}$ satisfies the equation
 \begin{equation}
 \frac{d}{dt}\Big(\frac{dT}{dP}\Big)=a\Big(\frac{dT}{dP}\Big)+b.
 \end{equation}
 The constants $a$ and $b$ can be estimated off line via simulation  and used to solve the latter equation.
 Moreover, if the settling time of this equation is  shorter than the control cycles
 then we  just use the steady-state value of Equation (16), which is $-\frac{b}{a}$. We feel confident that this additional
 approximation simplifies the control algorithm without significantly degrading its tracking performance.  Details of the computation of this term will be presented in the next section, where its effectiveness in
temperature control will be demonstrated.

 Using the above models for the systems $S_{1}$, $S_{2}$, and $S_{3}$, we can approximate the derivative term
 $\frac{dT}{d\phi}$ that is required by the regulation law via Equation  (5).
In fact,
combining  Equations (\ref{VoltageEquation}), (\ref{PowerEquation}), (\ref{StaticPowerEquation2}), and (\ref{DynamicPowerEquation}), and taking derivatives, we obtain, after some algebra, that
\begin{equation}
\frac{dT}{d\phi}=\frac{\Big(\frac{dT}{dP}\Big)\Big(m\frac{P_s}{V} + (P - P_s)\Big(\frac{1}{\phi}+\frac{2m}{V}\Big)\Big) }
{1-\Big(\frac{dT}{dP}\Big)P_s (\log 10) \Big(\frac{-\gamma}{T^2}\Big)}.
\end{equation}
We point out that all of the terms in the RHS of this equation except for
$P_s$ and $\frac{dT}{dP}$ can be obtained from real-time measurements of a core, $P_{s}$ can be calculated
online using Equation (\ref{StaticPowerEquation2}), and   $\frac{dT}{dP}$  can be estimated off-line by its steady-state
value, $-\frac{b}{a}$, obtained from (16).\\

Consider now the case of multiple cores on a chip, where the problem is to regulate their  temperatures
to given (not-necessary identical)  setpoints by adjusting their respective frequencies. Due to the thermal gradients between the cores,
it appears that their  temperatures have to be regulated jointly. However, extensive simulations,
described in the next section,
 revealed that the Jacobian matrix of the function relating the cores' frequency vector to the temperature vector is
diagonally dominant and this justifies the use of a distributed control where each core runs an adjustable-gain
integrator as described in Section II.
The details of this control law will be presented in the next section.

\section{SIMULATION EXPERIMENTS}

We tested the proposed controller on Manifold \cite{wang_ispass2014},
a cycle-level, full-system
processor simulation environment
with a suitable interface for injecting the thermal controller. The Manifold
framework simulates the architecture-level execution of applications
based on state-of-the-art physical models~\cite{song_ispass2014}. A
functional emulator front-end~\cite{kersey2012universal} boots a Linux
kernel and executes compiled binaries  from an established suite of
benchmarks~\cite{bienia2008parsec}.

The processor that we simulated consists of  four
out-of-order execution cores, a two-level cache hierarchy, and a memory
controller, and its architecture is shown in Figure 3. The centralized (joint) control consists of repeated
applications of  Equations
$(1)\rightarrow(8)\rightarrow(2)\rightarrow(3)$, where $u_{n-1}=\phi_{n-1}\in \mathds{R}^4$
is the vector of core frequencies during the $n^{th}$ cycle and $y_{n}=T_{n}\in \mathds{R}^4$ is the vector of core temperatures at the end of the $n^{th}$ cycle. Recall that Equation (8) denotes the controller's gain, and since it
is diagonal, the control is implemented by the cores in a distributed fashion. In contrast Equation (1) represents the
processor system and hence must be simulated jointly. This was done in Manifold in the following way.

\begin{figure}
\centering
\includegraphics[scale=0.8]{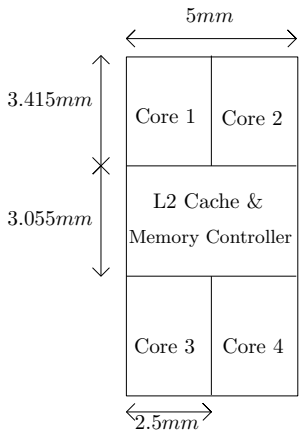}
\caption{Floor Plan of the 4 Core Processor}
\end{figure}

Equation (1) can be written as $T_{n}=g_{n}(\phi_{n-1})$, where $\phi_{n-1}:=(\phi_{n-1,1},\ldots,\phi_{n-1,4})^{\top}\in \mathds{R}^4$ and $T_{n}:=(T_{n,1},\ldots,T_{n,4})^{\top}\in \mathds{R}^4$ according to their respective co-ordinates, with
the second subscript $j=1,\ldots,4$ corresponding to the index of the core in Figure 3. In Equation (8) we approximate the $4\times 4$ Jacobian matrix $\frac{dT_{n}}{d\phi_{n-1}}$. Its diagonal terms, $\frac{\partial T_{n,j}}{\partial \phi_{n-1,j}}$,
$j=1,\ldots,4$, are just the terms $\frac{dT}{d\phi}$ in the Left-Hand Side (LHS) of Equation (17) with the
subscripts $n,j$ indicating  core $j$ at the $n^{th}$ control cycle.
As mentioned earlier all the terms in the RHS of (17) can be obtained from real-time
measurements and computation except for $\frac{dT}{dP}$, now referred to as
$\frac{dT_{n,j}}{dP_{n-1,j}}$. For estimating this term we used (16) in the steady state. To this end
we ran extensive Manifold simulations of the processor in open loop with various
input frequencies. Each simulation was run for successive cycles of $10$ms, long enough for the temperature
to reach its steady state, and it yielded traces of  power and its corresponding  temperature at  each cycle.
The traces, providing over $4,000$ data pairs per core,  indicated  a nearly-affine   power-to-temperature relation for each core
regardless of the physical state (frequencies and temperatures) at the other three cores. We  used the
MATLAB Curve-Fitting Toolbox to approximate these power-temperature  relations by respective lines,
whose slopes serve to estimate the terms $\frac{\partial T_{n,j}}{\partial P_{n-1,j}}$. Since the $P$-$T$ traces
were generated across the entire spectrum of   frequencies  at all four cores, the slopes of the approximating lines do not depend on $n$, although they may depend on $j=1,\ldots,4$ according to the processor's floor plan. Thus, the steady-state solution of Equation (16) in our case has the following approximation,
\begin{equation}
\frac{\partial T_{n,j}}{\partial P_{n-1,j}}\cong -\frac{b_{j}}{a_{j}},\ \ \ \ \ j=1,\ldots,4,
\end{equation}
whose right-hand side is the slope of the line associated with core $j$. The MATLAB Curve-Fitting Toolbox yielded the following values, $3.97, 5.242, 3.877, 4.055$ for cores $1 - 4$, respectively, with an R-Square confidence metric $> 0.97$. As a further
approximation we averaged these four numbers and thus used $-\frac{b_{j}}{a_{j}}\cong 4.286$ for $j=1,\ldots,4$.  This, in conjunction with (17) yields the terms $\frac{\partial T_{n,j}}{\partial\phi_{n-1,j}}$. We note that while this approximation of
$\frac{\partial T_{n,j}}{\partial P_{n-1,j}}$ is independent of $n$ or $j$, the partial derivative $\frac{\partial T_{n,j}}{\partial\phi_{n-1,j}}$ does depend on $n$ and $j$ through the other terms in the RHS of (17).

For the off-diagonal terms of $\frac{dT_{n}}{d\phi_{n-1}}$ we observe  (by the chain rule) that for $i,j=1,\ldots,4$,
\begin{equation}
\frac{\partial T_{n,i}}{\partial\phi_{n-1,j}}=\frac{\partial T_{n,i}}{\partial T_{n,j}}.\frac{\partial T_{n,j}}{\partial\phi_{n-1,j}}.
\end{equation}
The second multiplicative term in the RHS of (19) was discussed in the previous paragraph. As for the first term, we estimated it
by finite-difference approximations from the traces of simulation outputs. To this end we used HotSpot,
an established simulation platform designed to assess the thermal behavior of
digital designs \cite{hotspot2014}. The thermal model generated by HotSpot consists of  a linear, time-invariant circuit
comprised of
resistors and capacitors, where potentials and currents represent temperature and power, respectively.
The input to the circuit consists of current sources and the outputs are node voltages, and hence
HotSpot is a suitable tool for modeling the thermal behavior of the core.

Varying the input power to the cores one-at-a-time, we obtained the temperature variations from which the finite-difference
approximations for $\frac{\partial T_{n,i}}{\partial T_{n,j}}$ were derived. These approximating terms  also are independent of $n$
and hence denoted by $\frac{\partial T_{i}}{\partial T_{j}}$, but $\frac{\partial T_{n,j}}{\partial\phi_{n-1,j}}$ certainly depends on $n$ through the second term in the RHS of (19).\footnote{Manifold has the core frequencies as input but it does not permit us
to vary the core powers one-at-a-time, while HotSpot allows us to do just that. This is the reason
we used both simulation environments in the manner described above.}

The matrix $\frac{\partial T_{i}}{\partial T_{j}}$, $i,j=1,\ldots,4$, thus
obtained from HotSpot,   is
\begin{eqnarray}
\frac{\partial T_{i}}{\partial T_{j}} =
 \begin{bmatrix}
  1\times10^6 & 0.0439 & 0.003378 & 0.003378 \\
  0.0439 & 1\times10^6 & 0.003378 & 0.003378 \\
  0.003378  & 0.003378  & 1\times10^6 & 0.0439  \\
  0.003378 & 0.003378 & 0.0439 & 1\times10^6
 \end{bmatrix} \nonumber \\
 \times 10^{-6}. \nonumber
\end{eqnarray}
This is clearly diagonally dominant, and hence we expected the Jacobian matrix $\frac{dT_{n}}{d\phi_{n-1}}$ to be diagonally
dominant as well. This indeed was observed at each value of $n$, as the following randomly-chosen example
from our Manifold runs shows,
$$
\frac{dT_{n}}{d\phi_{n-1}} =
 \begin{bmatrix}
  23800 & 1109 & 73.78 & 72.73 \\
  1045 & 25270 & 73.78 & 72.73 \\
  80.405  & 85.37  & 21870 & 945  \\
  80.405 & 85.37 & 958 & 21530
 \end{bmatrix}
\times 10^{-6}.
$$
With this we felt confident in neglecting the off-diagonal terms of the Jacobian matrix, thereby replacing the joint core-temperature control
based on Equation (8) by four parallel one-dimensional controllers, one for each core, based on Equation (5).

\begin{figure*}[!t]
\centering
\includegraphics[width=\textwidth]{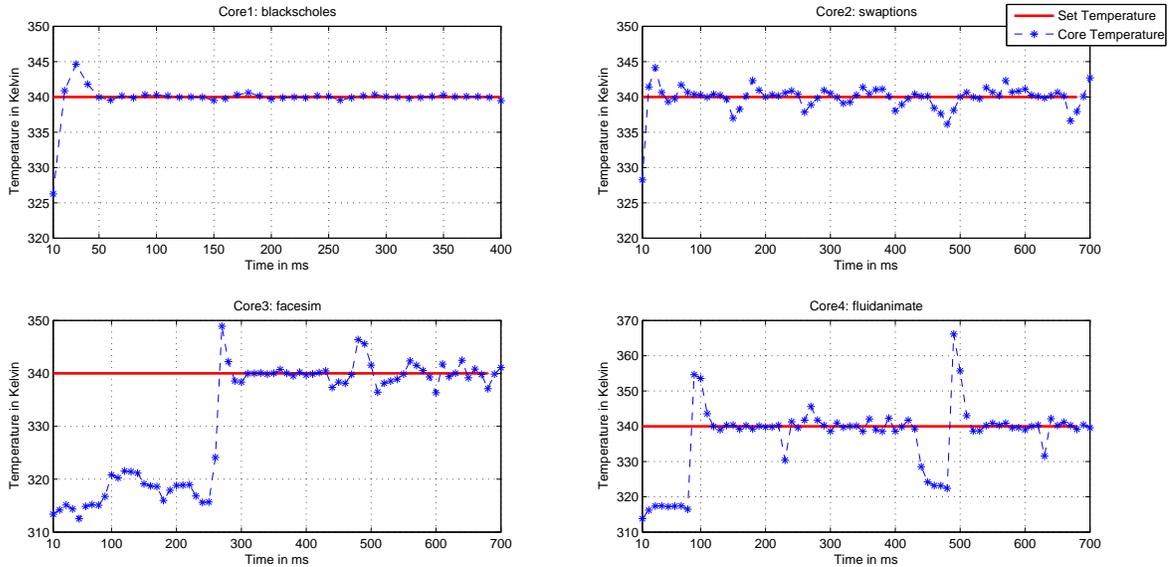}
\caption{Tracking results with Continuous Frequencies}
\end{figure*}

\begin{figure*}[!t]
\centering
\includegraphics[width=\textwidth]{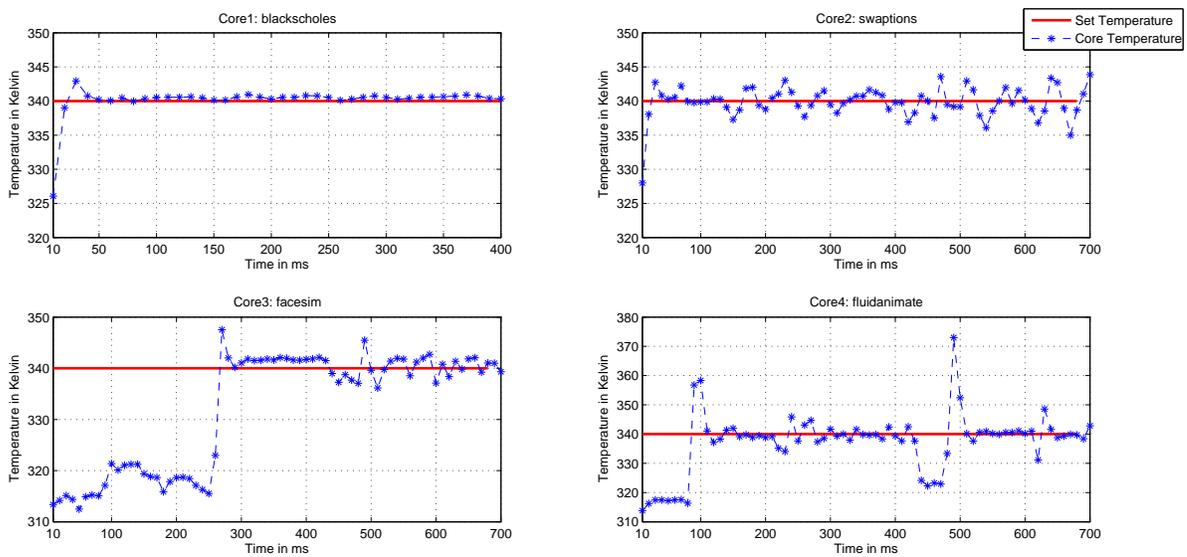}
\caption{Tracking results with Discrete Frequencies}
\end{figure*}

We implemented the distributed controller in conjunction with Manifold simulation of the processor. Each one of the cores
executed a different benchmark program from the parsec suite of benchmarks \cite{bienia2008parsec}: {\it blackscholes, swaptions, facesim, and fluidanimate} were executed   by Core $1$, Core $2$, Core $3$, and Core $4$ (see Figure 3), respectively.
The target temperature of all cores was set to $340$K, a typical  value, and the range of frequencies
was $1$GHz to $4.7$GHz. The control cycles at each one of the controllers were $10$ms.  {\it blackscholes} running on Core $1$ lasts $400$ms and hence the control was run for $40$ cycles, while the rest of the benchmarks take longer than $700$ms but we graph the results only for the first $70$ control cycles. The results are shown in the four graphs in Figure 4, and for each core we computed the average temperature from the end of the first overshoot to the cycle ending at the final time shown in the graph ($400$ ms for Core $1$, $700$ ms for the other cores).

In Core $1$ we notice convergence at $5$  iterations (control cycles) following a fast rise and a $5$-degree overshoot. The average temperature (from the end of the first overshoot to iteration $40$) is $339.995$K. In Core $2$ we see a similar rise and overshoot as in Core $1$, but then we note an oscillatory
behavior and not a smooth tracking. The reason is that the benchmark {\it swaptions} had large and rapid variations in its
activity factor $(\alpha(t))$ and hence in the dissipated dynamic power, causing ripples in the temperature profile. However, the computed average temperature is $339.96$K - arguably quite close to the target setpoint of 340K.

Core $3$  shows no tracking until $250$ms, then an overshoot followed by a $130$-ms smooth tracking, and a period of
minor ripples. The reason for the delayed tracking
is that during the first $250$ms the benchmark \textit{facesim} is in a data-fetch phase when  most of the computation units within the core are idle. Therefore there  is no significant dynamic power dissipation and the core temperature does not rise. During
that phase  the core frequency first climbs to its maximum  value ($4.7$Ghz) and then stays there
until time  $250$ms. Once the program enters the computation phase (time $> 250$ms), the dynamic power rises which causes the core temperature to increase and the controller is now able to track the set temperature of $340$K. The average temperature,
computed as before, was $340.204$K.

In Core $4$ the  benchmark program has two data-fetch periods and also periods of   wide-range  power dissipation during its execution. We discern a similar delayed tracking  as was observed with Core 3 but for a shorter duration, ending
at $t=80$ms.  Later the program enters another data-fetch phase in the time range of  $400$ - $500$ms,  causing the core temperature to drop while the frequency rises to its maximum value. In both cases the data-fetch  phase is followed by a computation phase which results in a temperature overshoot followed by a period of tracking except for ripples that are due to large variability in the dynamic power. The average temperature from the end of the first overshoot to the last  control cycle
 shown in the graph was  $339.565$K.

In the previous simulation we allowed the frequency to take any value in the range $1$GHz to $4.7$GHz. However, in a typical processor only a finite set of   frequencies  can be applied to a core. Therefore we repeated the simulation of the control technique for the following set of allowed frequencies, $\{1,1.5,1.8,3.4,3.7,3.9,4.0,4.1,4.2,4.4,4.7\}$ GHz. The only
 difference from the previous simulation is that in Equation (2) we took the control $u_{n}$ to be the
 nearest element in this set to the
computed term $u_{n-1}+A_{n}e_{n-1}$. The results are shown in Figure 5, and they are similar to those in Figure 4 except that slightly larger ripples and minor steady-state errors are discerned.
These  were expected, and are due to the quantization errors in the selection of frequencies.  However, the average temperatures  at the cores, from the end of the first overshoot to the final time, are quite close to the setpoint
reference:  $340.482$K, $339.986$K, $340.623$K,  and $339.392$K at Cores $1 - 4$, respectively.

We close this section by comparing the tracking performance of our adaptive-gain controller with those using fixed gains.
The need for an adaptive-gain control  arises from
unpredictable program activity factors $(\alpha(t))$, which may vary widely during the program.
We simulated the four-core system but applied the controllers only to core 4  running the
{\it fluidanimate} benchmark. The frequency range is continuous. We chose a low gain of
10 and a high gain of 120.
The graphs of the temperature traces obtained from these two gains as well as the variable-gain control  are shown in Figure 6. It is readily seen that the low gain results in the longest settling times, while
the high gain yields larger oscillations.  Not surprisingly, the tracking performance of the variable-gain controller is
better than those of the two fixed-gain controls.

\begin{figure}[!t]
\centering
\includegraphics[width=3.4in]{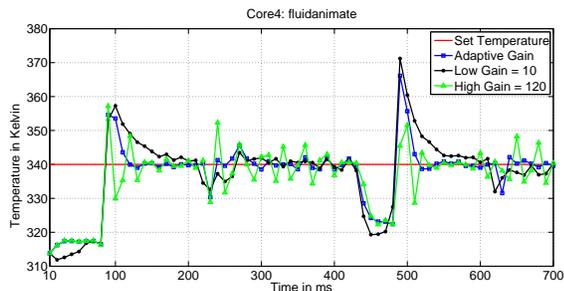}
\caption{Tracking results with fixed gains and variable gains}
\end{figure}

\section{CONCLUDING REMARKS}
\label{Section5}

Temperature regulation has emerged as a fundamental requirement of
modern and future processors. The state of the practice to date has
been dominated by ad-hoc adaptive heuristics. More recent attempts
have begun to apply the rich landscape of control theory to this
problem. However, these techniques have primarily dealt with
temperature as a constraint while controlling power dissipation.

This paper makes a subtle but important observation - temperature
ought to be directly regulated to track a target value while power
should be managed to maximize performance. Regulating chip-wide
temperature to a balanced  thermal field is necessary
while preventing transitions across a maximum temperature, since  the latter can
produce thermal fields that adversely affect reliability and
performance. Furthermore, unlike prior works we consider a  nonlinear,
time-varying plant model that explicitly captures the exponential
dependence of temperature and static power, and devise a distributed control technique that trades off precision with simplicity of
real-time computations.   Simulation results using a
full system, cycle level simulator executing industry standard
benchmark programs indicate  convergence of our regulation
technique despite the modeling approximations.


\begin{thebibliography}{10}
\providecommand{\url}[1]{#1}
\csname url@samestyle\endcsname
\providecommand{\newblock}{\relax}
\providecommand{\bibinfo}[2]{#2}
\providecommand{\BIBentrySTDinterwordspacing}{\spaceskip=0pt\relax}
\providecommand{\BIBentryALTinterwordstretchfactor}{4}
\providecommand{\BIBentryALTinterwordspacing}{\spaceskip=\fontdimen2\font plus
\BIBentryALTinterwordstretchfactor\fontdimen3\font minus
  \fontdimen4\font\relax}
\providecommand{\BIBforeignlanguage}[2]{{%
\expandafter\ifx\csname l@#1\endcsname\relax
\typeout{** WARNING: IEEEtran.bst: No hyphenation pattern has been}%
\typeout{** loaded for the language `#1'. Using the pattern for}%
\typeout{** the default language instead.}%
\else
\language=\csname l@#1\endcsname
\fi
#2}}
\providecommand{\BIBdecl}{\relax}
\BIBdecl

\bibitem{dennard1974design}
R.~H. Dennard, F.~H. Gaensslen, V.~L. Rideout, E.~Bassous, and A.~R. LeBlanc,
  ``Design of ion-implanted mosfet's with very small physical dimensions,''
  \emph{Solid-State Circuits, IEEE Journal of}, vol.~9, no.~5, pp. 256--268,
  1974.

\bibitem{ITRS2011}
``International Technology Roadmap for Semiconductors (ITRS), 2011,''
  \url{http://www.itrs.net/Links/2011ITRS/Home2011.htm}, accessed: 2014-03-15.

\bibitem{skadron2003temperature}
K.~Skadron, M.~R. Stan, W.~Huang, S.~Velusamy, K.~Sankaranarayanan, and
  D.~Tarjan, ``Temperature-aware microarchitecture,'' in \emph{ACM SIGARCH
  Computer Architecture News}, vol.~31, no.~2.\hskip 1em plus 0.5em minus
  0.4em\relax ACM, 2003, pp. 2--13.

\bibitem{skadron2002control}
K.~Skadron, T.~Abdelzaher, and M.~R. Stan, ``Control-theoretic techniques and
  thermal-rc modeling for accurate and localized dynamic thermal management,''
  in \emph{High-Performance Computer Architecture, 2002. Proceedings. Eighth
  International Symposium on}.\hskip 1em plus 0.5em minus 0.4em\relax IEEE,
  2002, pp. 17--28.

  \bibitem{liu2012neighbor}
G.~Liu, M.~Fan, and G.~Quan, ``Neighbor-aware dynamic thermal management for
  multi-core platform,'' in \emph{Design, Automation \& Test in Europe
  Conference \& Exhibition (DATE), 2012}.\hskip 1em plus 0.5em minus
  0.4em\relax IEEE, 2012, pp. 187--192.

\bibitem{yeo2008predictive}
I.~Yeo, C.~C. Liu, and E.~J. Kim, ``Predictive dynamic thermal management for
  multicore systems,'' in \emph{Proceedings of the 45th annual Design
  Automation Conference}.\hskip 1em plus 0.5em minus 0.4em\relax ACM, 2008, pp.
  734--739.

  \bibitem{kong2012recent}
J.~Kong, S.~W. Chung, and K.~Skadron, ``Recent thermal management techniques
  for microprocessors,'' \emph{ACM Computing Surveys (CSUR)}, vol.~44, no.~3,
  p.~13, 2012.

\bibitem{donald2006techniques}
J.~Donald and M.~Martonosi, ``Techniques for multicore thermal management:
  Classification and new exploration,'' \emph{ACM SIGARCH Computer Architecture
  News}, vol.~34, no.~2, pp. 78--88, 2006.

\bibitem{qian2011cyber}
H.~Qian, X.~Huang, H.~Yu, and C.~H. Chang, ``Cyber-physical thermal management
  of 3d multi-core cache-processor system with microfluidic cooling,''
  \emph{Journal of Low Power Electronics}, vol.~7, no.~1, pp. 110--121, 2011.

\bibitem{shi2010dynamic}
B.~Shi, Y.~Zhang, and A.~Srivastava, ``Dynamic thermal management for single
  and multicore processors under soft thermal constraints,'' in
  \emph{Proceedings of the 16th ACM/IEEE international symposium on Low power
  electronics and design}.\hskip 1em plus 0.5em minus 0.4em\relax ACM, 2010,
  pp. 165--170.

\bibitem{fu2012feedback}
Y.~Fu, N.~Kottenstette, C.~Lu, and X.~D. Koutsoukos, ``Feedback thermal control
  of real-time systems on multicore processors,'' in \emph{Proceedings of the
  tenth ACM international conference on Embedded software}.\hskip 1em plus
  0.5em minus 0.4em\relax ACM, 2012, pp. 113--122.

\bibitem{zanini2009multicore}
F.~Zanini, D.~Atienza, L.~Benini, and G.~De~Micheli, ``Multicore thermal
  management with model predictive control,'' in \emph{Circuit Theory and
  Design, 2009. ECCTD 2009. European Conference on}.\hskip 1em plus 0.5em minus
  0.4em\relax IEEE, 2009, pp. 711--714.

\bibitem{wang2011adaptive}
X.~Wang, K.~Ma, and Y.~Wang, ``Adaptive power control with online model
  estimation for chip multiprocessors,'' \emph{Parallel and Distributed
  Systems, IEEE Transactions on}, vol.~22, no.~10, pp. 1681--1696, 2011.

\bibitem{murali2008temperature}
S.~Murali, A.~Mutapcic, D.~Atienza, R.~Gupta, S.~Boyd, L.~Benini, and
  G.~De~Micheli, ``Temperature control of high-performance multi-core platforms
  using convex optimization,'' in \emph{Design, Automation and Test in Europe,
  2008. DATE'08}.\hskip 1em plus 0.5em minus 0.4em\relax IEEE, 2008, pp.
  110--115.

\bibitem{bartolini2013thermal}
A.~Bartolini, M.~Cacciari, A.~Tilli, and L.~Benini, ``Thermal and energy
  management of high-performance multicores: Distributed and self-calibrating
  model-predictive controller,'' \emph{Parallel and Distributed Systems, IEEE
  Transactions on}, vol.~24, no.~1, pp. 170--183, 2013.



\bibitem{rotem2012power}
E.~Rotem, A.~Naveh, D.~Rajwan, A.~Ananthakrishnan, and E.~Weissmann,
  ``Power-management architecture of the intel microarchitecture code-named
  sandy bridge,'' \emph{IEEE Micro}, pp. 20--27, 2012.

\bibitem{AMDPhenom}
``AMD Phenom II Key Architectural Features,''
  \url{http://www.amd.com/us/products/desktop/processors/phenom-ii/Pages/phenom-ii-key-architectural-features.aspx},
  accessed: 2014-03-15.


\bibitem{paul2013}
I.~Paul, S.~Manne, L.~Bircher, M.~Arora, and S.~Yalamanchili, “Cooperative Boosting: Needy vs. Greedy Power Management,” \emph{IEEE/ACM International Symposium on Computer Architecture (ISCA)}, June 2013.

\bibitem{song2015}
W.~Song, S.~Mukhopadhyay, and S.~Yalamanchili, “Architectural Reliability: Lifetime Reliability Characterization and Management for Many Core Processors, “ \emph{IEEE Computer Architecture Letters}, to appear.



\bibitem{almoosa2012power}
N.~Almoosa, W.~Song, Y.~Wardi, and S.~Yalamanchili, ``A power capping
  controller for multicore processors,'' in \emph{American Control Conference
  (ACC), 2012}.\hskip 1em plus 0.5em minus 0.4em\relax IEEE, 2012, pp.
  4709--4714.


\bibitem{lancaster1966error}
P.~Lancaster ``Error analysis for the Newton-Raphson method,'' in \emph{Numerische Mathematik, 1966},
vol.~9, pp. 55--68, 1966.

\bibitem{ortega2000iterative}
\BIBentryALTinterwordspacing
J.~Ortega, and Rheinboldt, C.~Werner \emph{Iterative solution of nonlinear equations in several variables},
Siam, 2000.
\BIBentrySTDinterwordspacing

\bibitem{burd2000dynamic}
T.~D. Burd, T.~A. Pering, A.~J. Stratakos, and R.~W. Brodersen, ``A dynamic
  voltage scaled microprocessor system,'' \emph{Solid-State Circuits, IEEE
  Journal of}, vol.~35, no.~11, pp. 1571--1580, 2000.

\bibitem{mcgowen2006power}
R.~McGowen, C.~A. Poirier, C.~Bostak, J.~Ignowski, M.~Millican, W.~H. Parks,
  and S.~Naffziger, ``Power and temperature control on a 90-nm itanium family
  processor,'' \emph{Solid-State Circuits, IEEE Journal of}, vol.~41, no.~1,
  pp. 229--237, 2006.

\bibitem{butts2000static}
J.~A. Butts and G.~S. Sohi, ``A static power model for architects,'' in
  \emph{Proceedings of the 33rd annual ACM/IEEE international symposium on
  Microarchitecture}.\hskip 1em plus 0.5em minus 0.4em\relax ACM, 2000, pp.
  191--201.

\bibitem{rabaey2009low}
\BIBentryALTinterwordspacing
J.~Rabaey, \emph{Low Power Design Essentials}, ser. Integrated Circuits and
  Systems.\hskip 1em plus 0.5em minus 0.4em\relax Springer, 2009. [Online].
  Available: \url{http://books.google.com/books?id=A-sBy\_nmQ8wC}
\BIBentrySTDinterwordspacing

\bibitem{han2007tilts}
Y.~Han, I.~Koren, and C.~M. Krishna, ``Tilts: A fast architectural-level
  transient thermal simulation method,'' \emph{Journal of Low Power
  Electronics}, vol.~3, no.~1, pp. 13--21, 2007.

\bibitem{wang_ispass2014} Wang et al., ``Manifold: A Parallel Simulation Framework for Multicore Systems,'' \textit{ISPASS}, Mar. 2014.

\bibitem{song_ispass2014} Song et al., ``Energy Introspector: A Parallel, Composable Framework for Integrated Power-Reliability-Thermal Modeling for Multicore Architectures,'' \textit{ISPASS}, Mar. 2014.

\bibitem{kersey2012universal}
C.~D. Kersey, A.~Rodrigues, and S.~Yalamanchili, ``A universal parallel
  front-end for execution driven microarchitecture simulation,'' in
  \emph{Proceedings of the 2012 Workshop on Rapid Simulation and Performance
  Evaluation: Methods and Tools}.\hskip 1em plus 0.5em minus 0.4em\relax ACM,
  2012, pp. 25--32.

\bibitem{bienia2008parsec}
C.~Bienia, S.~Kumar, and K.~Li, ``Parsec vs. splash-2: A quantitative
  comparison of two multithreaded benchmark suites on chip-multiprocessors,''
  in \emph{Workload Characterization, 2008. IISWC 2008. IEEE International
  Symposium on}.\hskip 1em plus 0.5em minus 0.4em\relax IEEE, 2008, pp. 47--56.

  \bibitem{hotspot2014}
  ``HotSpot Version 5.0''
    \url{http://lava.cs.virginia.edu/HotSpot/index.htm}, accessed: 2014-09-19.






\end{thebibliography}
\end{document}